\documentclass{article}
\usepackage{amsfonts}
\usepackage{amssymb}
\usepackage{amsmath}

\newtheorem{theorem}{Theorem}

\newtheorem{corollary}[theorem]{Corollary}

\newtheorem{definition}[theorem]{Definition}

\newtheorem{lemma}[theorem]{Lemma}

\newtheorem{proposition}[theorem]{Proposition}
\newtheorem{remark}[theorem]{Remark}

\newenvironment{proof}[1][Proof]{\noindent\textbf{#1.} }{\ \rule{0.5em}{0.5em}}
\input{tcilatex}

\begin{document}

\title{The Burkholder-Davis-Gundy Inequality for Enhanced Martingales}
\author{Peter Friz \and Nicolas Victoir}
\maketitle

\begin{abstract}
Multi-dimensional continuous local martingales, enhanced with their
stochastic area process, give rise to geometric rough paths with a.s. finite
homogenous $p$-variation, $p>2$. Here we go one step further and establish
quantitative bounds of the $p$-variation norm in the form of a BDG
inequality. Our proofs are based on old ideas by L\'{e}pingle. We also
discuss geodesic and piecewise linear approximations.
\end{abstract}

\section{\protect\bigskip Introduction}

The theory of rough paths provides a new and robust way to drive
differential equations by multi-dimensional stochastic processes in a 
\textit{deterministic way}. In most cases, this is achieved by taking into
account a certain stochastic area process and by establishing fine
regularity properties of the resulting \textit{enhanced} process. The object
of study in this paper is a $d$-dimensional continuous local martingale $M$
null at $0$ for which the area is defined by iterated stochastic
integration; the area process $A_{t}$ is simply the anti-symmetric part of
the iterated Stratonovich integral,%
\begin{equation*}
\mathbf{M}_{t}^{2}\equiv \int_{0}^{t}\int_{0}^{s}dM_{r}\otimes \circ
dM_{s}\in \mathbb{R}^{d}\otimes \mathbb{R}^{d}.
\end{equation*}%
Note that the symmetric part of $\mathbf{M}_{t}^{2}$ is given by $\frac{1}{2}%
M_{t}\otimes M_{t}$ and hence redundant if one knows $\mathbf{M}%
_{t}^{1}\equiv M_{t}$. It follows that the enhanced process $\mathbf{M}%
\equiv \left( 1,\mathbf{M}^{1},\mathbf{M}^{2}\right) \in \mathbb{R\oplus }%
\mathbb{R}^{d}\oplus \mathbb{R}^{d}\otimes \mathbb{R}^{d}$ lives in
submanifold, namely in $G^{2}\left( \mathbb{R}^{d}\right) \equiv \exp \left( 
\mathbb{R}^{d}\oplus so\left( d\right) \right) $, where $\exp :\left(
x,a\right) \mapsto \left( 1,x,a+\frac{1}{2}x\otimes x\right) .$ The space $%
\mathbb{R}^{d}\oplus so\left( d\right) $ carries a Lie algebra structure and
induces a (Lie-)group structure on $G^{2}\left( \mathbb{R}^{d}\right) $. The
interest in this algebraic exercise is that the resulting product operation
on $G^{2}\left( \mathbb{R}^{d}\right) $ is exactly what one needs to patch
together "iterated integral increments" over adjacent intervals. $%
G^{2}\left( \mathbb{R}^{d}\right) $ is also a metric (in fact,\ Polish)
space under the Carnot-Caratheodory metric $d$. Intuitively, the distance of
two points under this metric is the length of the shortest path in $\mathbb{R%
}^{d}$ which wipes out a prescribed area. When $d=2$, geodesics are seen to
be parts of circles. $G^{2}\left( \mathbb{R}^{d}\right) $ carries a dilation
induced by $\left( x,a\right) \mapsto \left( \lambda x,\lambda ^{2}a\right) $
for real $\lambda $. In fact, the CC-metric is induced by a sub-additive
norm, homogenuous w.r.t. dilation. Since all continuous homogenous norms are
Lipschitz equivalent, computations are often carried out w.r.t. $\left\vert
\left\vert \left\vert \left( x,a\right) \right\vert \right\vert \right\vert =
$ $\left\vert x\right\vert +\left\vert a\right\vert ^{1/2}$. We refer to 
\cite{Ly98, LQ02} for background on rough paths, \cite{FV1} contains a more
detailed discussion of the relevant geometry and algebra. The notion of a
(weak) geometric $p$-rough path \cite{FV2} becomes quite elegant: by
definition, one requires that the $G^{2}\left( \mathbb{R}^{d}\right) $%
-valued path $\mathbf{M}$ has finite $p$-variation%
\begin{equation*}
\left\Vert \mathbf{M}\right\Vert _{p\text{-var;}[0,T]}=\left( \sup_{0\leq
t_{1}\leq ...\leq t_{n}\leq T}\sum d\left( \mathbf{M}_{t_{i}}\mathbf{,M}%
_{t_{i+1}}\right) ^{p}\right) ^{1/p}<\infty .
\end{equation*}%
Is is known\ \cite{CL} that this holds for a.e. $\mathbf{M}=\mathbf{M}\left(
\omega \right) $ when $p>2$. The first main topic of this paper is to
establish quantitative bounds of the $p$-variation norm in the form of a
two-sided BDG inequality: for any moderate function $F$ such as $x\mapsto
x^{r}$ for $r>0$,%
\begin{equation*}
\mathbb{E}\left( F\left( \left\Vert \mathbf{M}\right\Vert _{p\text{-var;}%
[0,T}\right) \right) \sim \mathbb{E}\left( F\left( \left\vert \left\langle
M\right\rangle _{T}\right\vert ^{1/2}\right) \right) .
\end{equation*}%
The algebraic and geometric preparations made above prove crucial to recycle
many of the arguments given in L\'{e}pingle's seminal paper \cite{L} from
1976. Secondly, we discuss approximations and show $L^{q}$-convergence (at
least for $q>1$) of lifted piecewise linear approximations of a continuous $%
L^{q}$-martingale w.r.t. homogenous $p$-variation topology.

The authors would like to thank D. L\'{e}pingle for a helpful email
exchange. 

\section{Preliminaries}

We write $\mathcal{M}_{0,\text{\textrm{loc}}}^{c}\left( [0,\infty ),\mathbb{R%
}^{d}\right) $ or $\mathcal{M}_{0,\text{\textrm{loc}}}^{c}\left( \mathbb{R}%
^{d}\right) $ for the class of $\mathbb{R}^{d}$-valued continuous local
martingales $M:[0,\infty )\rightarrow \mathbb{R}^{d}$ null at $0$. The
bracket process $\left\langle M\right\rangle :[0,\infty )\rightarrow \mathbb{%
R}^{d}$ is defined component-wise, the $i^{th}$ component is given by the
usual bracket $\left\langle M^{i}\right\rangle =\left\langle
M^{i},M^{i}\right\rangle .$

The area-process $A:[0,\infty )\rightarrow so\left( d\right) $ is defined by
It\^{o}- or Stratonovich stochastic integration. As the matrix $\left\langle
M^{i},M^{j}\right\rangle $ is symmetric both lead to the \textit{same} area,%
\begin{eqnarray*}
A_{t}^{i,j} &=&\frac{1}{2}\left(
\int_{0}^{t}M^{i}dM^{j}-\int_{0}^{t}M^{j}dM^{i}\right) \\
&=&\frac{1}{2}\left( \int_{0}^{t}M^{i}\circ dM^{j}-\int_{0}^{t}M^{j}\circ
dM^{i}\right) .
\end{eqnarray*}%
We note that the area-process is a vector-valued continuous martingale. By
disregarding a null-set we can and will assume that $M$ and $A$ are
continuous.

\begin{definition}
Set $S_{2}\left( M\right) :=\mathbf{M}:=\exp \left( M+A\right) $ so that $%
\mathbf{M}\in C\left( [0,\infty ),G^{2}\left( \mathbb{R}^{d}\right) \right) $%
. The resulting class of enhanced (continuous, local) martingales is denoted
by $\mathcal{M}_{0,\text{\textrm{loc}}}^{c}\left( G^{2}\left( \mathbb{R}%
^{d}\right) \right) $. We refer to the operation $S_{2}:M\mapsto \mathbf{M}$
as lift.
\end{definition}

The lift is compatible with the stopping and time-changes.

\begin{lemma}
\bigskip (i) Let $\tau $ be a stopping time. Then $\mathbf{M}^{\tau
}=S_{2}\left( M^{\tau }\right) $. (ii) Let $\phi $ be a time-change, that
is, a family $\phi _{s},\,s\geq 0,$ of stopping times such that the maps $%
s\mapsto \phi _{s}$ are a.s. increasing and right-continuous. Assume that $M$
is constant on each interval $\left[ \phi _{t-},\phi _{t}\right] $. Then $%
M\circ \phi $ is a continuous local martingale and%
\begin{equation*}
\mathbf{M}\circ \phi =S_{2}\left( M\circ \phi \right) .
\end{equation*}
\end{lemma}

\begin{proof}
Stopped processes are special cases of time-changed processes (take $\phi
_{t}=t\wedge \tau $) so it suffices to show the second statement. To this
end, recall the compatibility of a time change $\phi $ and stochastic
integration w.r.t. a continuous local martingale, constant on each interval $%
\left[ \phi _{t-},\phi _{t}\right] $, Proposition V.1.5. (ii)\ of \cite{RY}.
The lift is a special case of stochastic integration.
\end{proof}

The lift is also compatible with respect to scaling and concatenation of
(local martingale) paths.

\begin{lemma}
(i) If $\delta _{c}:$ $G^{2}\left( \mathbb{R}^{d}\right) \rightarrow
G^{2}\left( \mathbb{R}^{d}\right) $ is the dilation operator given by $%
\delta _{c}\exp \left( x+a\right) =\exp \left( cx+c^{2}a\right) $ then $%
\delta _{c}\mathbf{M}=S_{2}\left( cM\right) $. (ii) We have%
\begin{equation*}
S_{2}\left( M\right) _{0,t}=S_{2}\left( \left. M\right\vert _{\left[ 0,s%
\right] }\ast \left. M\right\vert _{\left[ s,t\right] }\right)
_{0,t}=S_{2}\left( M\right) _{0,s}\otimes S_{2}\left( M\right)
_{s,t},\,\,\,0\leq s\leq t<\infty
\end{equation*}
\end{lemma}

\begin{proof}
(i) follows is trivial consequence of linearity of stochastic integrals.
(ii) is true whenever a first order calculus underlies the lift. It now
suffices to note that $S_{2}$ is (equivalently) defined as Stratonovich
lift, 
\begin{equation*}
S_{2}\left( M\right) _{t}=\exp \left( M_{t}+A_{t}\right)
=1+M_{t}+\int_{0}^{t}M_{s}\otimes \circ dM_{s}\text{.}
\end{equation*}
\end{proof}

\begin{definition}
$F:\mathbb{R}^{+}\rightarrow \mathbb{R}^{+}$ is \textit{moderate} if (i) $F$
is continuous and increasing, (ii) $F\left( x\right) =0$ if and only if $x=0$
and (iii) for some (and then for every) $\alpha >1$,%
\begin{equation*}
\sup_{x>0}\frac{F\left( \alpha x\right) }{F\left( x\right) }<\infty .
\end{equation*}
\end{definition}

The following result is found, for instance, in \cite[p93]{RW2}.

\begin{theorem}[Burkholder-Davis-Gundy ]
Let $F$ be a moderate function, $M\in \mathcal{M}_{0,\text{\textrm{loc}}%
}^{c}\left( \mathbb{R}\right) $. Then there exists a constant $C=C\left(
F,d,\left\vert \cdot \right\vert \right) $ so that 
\begin{equation*}
C^{-1}\mathbb{E}\left( F\left( \left\vert \left\langle M\right\rangle
_{\infty }\right\vert ^{1/2}\right) \right) \leq \mathbb{E}\left( F\left(
\sup_{s\geq 0}\left\vert M_{s}\right\vert \right) \right) \leq C\mathbb{E}%
\left( F\left( \left\vert \left\langle M\right\rangle _{\infty }\right\vert
^{1/2}\right) \right) .
\end{equation*}
\end{theorem}

We collect a few properties of moderate functions.

\begin{lemma}
(i) $x\mapsto F\left( x\right) $ moderate iff $\mapsto F\left(
x^{1/2}\right) $ moderate.\newline
(ii) Given $c,A,B>0$ $:c^{-1}A\leq B\leq cA\implies \exists C=C\left(
c,F\right) :$%
\begin{equation*}
C^{-1}F\left( A\right) \leq F\left( B\right) \leq CF\left( A\right) .
\end{equation*}%
\newline
(iii) $\exists C:\forall x,y>0:F\left( x+y\right) \leq C\left[ F\left(
x\right) +F\left( y\right) \right] .$
\end{lemma}

\begin{proof}
(i),(ii) are left to the reader. Ad (iii): W.l.o.g. \thinspace $x<y$, then $%
F\left( x+y\right) \leq F\left( 2y\right) \leq CF\left( y\right) $ by
moderate growth of $F$.
\end{proof}

\begin{corollary}
Let $F$ be a moderate function, $M\in \mathcal{M}_{0,\text{\textrm{loc}}%
}^{c}\left( \mathbb{R}^{d}\right) $ and $\left\vert \cdot \right\vert $ a
norm on $\mathbb{R}^{d}$. Then there exists a constant $C=C\left(
F,d,\left\vert \cdot \right\vert \right) $ so that 
\begin{equation*}
C^{-1}\mathbb{E}\left( F\left( \left\vert \left\langle M\right\rangle
_{\infty }\right\vert ^{1/2}\right) \right) \leq \mathbb{E}\left( F\left(
\sup_{s\geq 0}\left\vert M_{s}\right\vert \right) \right) \leq C\mathbb{E}%
\left( F\left( \left\vert \left\langle M\right\rangle _{\infty }\right\vert
^{1/2}\right) \right) .
\end{equation*}
\end{corollary}

\begin{proof}
When $\left\vert a\right\vert =\max \left\{ \left\vert a^{1}\right\vert
,...,\left\vert a^{d}\right\vert \right\} $ this is a simple consequence of
BDG for $\mathcal{M}_{0,\text{\textrm{loc}}}^{c}\left( \mathbb{R}\right) $,
applied componentwise. The lemma above shows that one can switch to
Lipschitz equivalent norms.
\end{proof}

From L\'{e}pingle \cite{L}, $\sup_{s\geq 0}\left\vert M_{s}\right\vert $
above can be replaced by the $p$-variation norm\footnote{%
In the next section, we will see a more general version of this.}. Noting
that the $p$-variation of a discrete-time martingale $\left( Y_{n}\right) $
is naturally defined as%
\begin{equation*}
\left\vert Y\right\vert _{p\text{-var}}\equiv \left[ \sup_{\left(
n_{k}\right) \nearrow }\sum_{k}\left\vert Y_{n_{k+1}}-Y_{n_{k}}\right\vert
^{p}\right] ^{1/p},
\end{equation*}%
the following lemma is best viewed as a BDG-type upper bound for
discrete-time martingales.

\begin{lemma}
\label{Prop2bLep76} Let $F$ be moderate. If $1<q<p\leq 2$ or $1=q=p$ then
there exists a constant $c$ such that for all, possibly $\mathbb{R}^{d}$%
-valued, discrete-time martingales $\left( Y_{n}:n\in \mathbb{Z}^{+}\right) $%
\begin{equation*}
\mathbb{E}\left( F\left( \left\vert Y\right\vert _{p\text{-var}}\right)
\right) \leq c\mathbb{E}\left[ F\left( \left[ \sum_{n}\left\vert
Y_{n+1}-Y_{n}\right\vert ^{q}\right] ^{1/q}\right) \right] .
\end{equation*}
\end{lemma}

\begin{proof}
For $d=1$ we can use Proposition 2.b in \cite{L} with the remark that a
discrete-time martingale can be viewed as a particular case of a
continuous-time martingale with purely discontinuous sample paths. As above,
the extension to $d>1$ does not pose any difficulty.
\end{proof}

\section{\protect\bigskip BDG on the group}

\begin{lemma}[{Good$\,\protect\lambda \,$inequality, \protect\cite[p.94]{RW2}%
}]
Let $X,Y$ be nonnegative random variables, and suppose there exists $\beta
>1 $ such that for all $\lambda >0,\delta >0,$%
\begin{equation*}
\mathbb{P}\left( X>\beta \lambda ,Y<\delta \lambda \right) \leq \psi \left(
\delta \right) \mathbb{P}\left( X>\lambda \right)
\end{equation*}%
where $\psi \left( \delta \right) \searrow 0$ when $\delta \searrow 0$.
There, for each moderate function $F,$ there exists a constant $C$ depending
only on $\beta ,\psi ,F$ such that%
\begin{equation*}
\mathbb{E}\left( F\left( X\right) \right) \leq C\mathbb{E}\left( F\left(
Y\right) \right) .
\end{equation*}
\end{lemma}

\begin{proposition}
\label{TschebTypeBound}Let $\left\vert \cdot \right\vert ,\,\left\Vert \cdot
\right\Vert $ continuous homogonous norm on $\mathbb{R}^{d},G^{2}\left( 
\mathbb{R}^{d}\right) $ respectively. Then there exists a constant $%
A=A\left( d,\left\vert \cdot \right\vert ,\left\Vert \cdot \right\Vert
\right) $ such that 
\begin{equation*}
\forall \mathbf{M}\in \mathcal{M}_{0,\text{\textrm{loc}}}^{c}\left(
G^{2}\left( \mathbb{R}^{d}\right) \right) \forall \lambda >0:\mathbb{P}%
\left( \sup_{u,v\geq 0}\left\Vert \mathbf{M}_{u,v}\right\Vert \geq \lambda
\right) \leq A\frac{\mathbb{E}\left( \left\vert \left\langle M\right\rangle
_{\infty }\right\vert \right) }{\lambda ^{2}}.
\end{equation*}
\end{proposition}

\begin{proof}
We note that $\sup_{u,v\geq 0}\left\Vert \mathbf{M}_{u,v}\right\Vert \leq
2\sup_{t\geq 0}\left\Vert \mathbf{M}_{t}\right\Vert $. By equivalence of
homogeneous norm, 
\begin{equation*}
\left\Vert \mathbf{M}_{t}\right\Vert ^{2}\leq C\left( \left\vert
M_{t}\right\vert ^{2}+\left\vert A_{t}\right\vert \right) .
\end{equation*}%
From BDG, $\mathbb{E}\left( \sup_{u\geq 0}\left\vert M_{u}\right\vert
^{2}\right) \leq C\mathbb{E}\left( \left\vert \left\langle M\right\rangle
_{\infty }\right\vert \right) .$ Note that $u\mapsto \left\vert \left\langle
M\right\rangle \right\vert _{u}:=\sum_{i=1}^{d}\left\langle
M^{i}\right\rangle _{u}$ is increasing (in fact, there is no loss in
generality in assuming that $\left\vert \cdot \right\vert $ on $\mathbb{R}%
^{d}$ is given given by $\left\vert a\right\vert =\sum \left\vert
a^{i}\right\vert $ ...). Then, using BDG again, 
\begin{eqnarray*}
\mathbb{E}\left( \sup_{u\geq 0}\left\vert A_{u}\right\vert \right) &\leq &C%
\mathbb{E}\left( \left\vert \int_{0}^{\infty }\left\vert M_{u}\right\vert
^{2}d\left\vert \left\langle M\right\rangle \right\vert _{u}\right\vert
^{1/2}\right) \\
&\leq &C\mathbb{E}\left( \sup_{u\geq 0}\left\vert M_{u}\right\vert
.\left\vert \left\langle M\right\rangle \right\vert _{\infty }^{1/2}\right)
\\
&\leq &C\sqrt{\mathbb{E}\sup_{u\geq 0}\left\vert M_{u}\right\vert ^{2}}\sqrt{%
\mathbb{E}\left[ \left\vert \left\langle M\right\rangle \right\vert _{\infty
}\right] } \\
&\leq &C\mathbb{E}\left( \left\vert \left\langle M\right\rangle \right\vert
_{\infty }\right) .
\end{eqnarray*}%
An application of Chebyshev's inequality finishes the proof.
\end{proof}

\begin{theorem}
\label{BDGgroupUniform}Let $F$ be a moderate function, $\mathbf{M}\in 
\mathcal{M}_{0,\text{\textrm{loc}}}^{c}\left( G^{2}\left( \mathbb{R}%
^{d}\right) \right) ,$ and $\left\vert \cdot \right\vert ,\,\left\Vert \cdot
\right\Vert $ continuous homogonous norm on $\mathbb{R}^{d},G^{2}\left( 
\mathbb{R}^{d}\right) $ respectively. Then there exists a constant $%
C=C\left( F,d,\left\vert \cdot \right\vert ,\left\Vert \cdot \right\Vert
\right) $ so that%
\begin{equation*}
C^{-1}\mathbb{E}\left( F\left( \left\vert \left\langle M\right\rangle
_{\infty }\right\vert ^{1/2}\right) \right) \leq \mathbb{E}\left( F\left(
\sup_{s,t\geq 0}\left\Vert \mathbf{M}_{s,t}\right\Vert \right) \right) \leq C%
\mathbb{E}\left( F\left( \left\vert \left\langle M\right\rangle _{\infty
}\right\vert ^{1/2}\right) \right) .
\end{equation*}
\end{theorem}

\begin{proof}
The lower bound comes from%
\begin{equation*}
\left\Vert \mathbf{M}_{s,t}\right\Vert \geq \left\vert M_{s,t}\right\vert
\end{equation*}%
the monotonicity of $F$ and the classical BDG lower bound. We prove the
upper-bound: we fix $\lambda ,\delta >0$ and $\beta >1,$ and we define the
stopping times 
\begin{eqnarray*}
S_{1} &=&\inf \left\{ t>0,\sup_{u,v\in \left[ 0,t\right] }\left\Vert \mathbf{%
M}_{u,v}\right\Vert >\beta \lambda \right\} , \\
S_{2} &=&\inf \left\{ t>0,\sup_{u,v\in \left[ 0,t\right] }\left\Vert \mathbf{%
M}_{u,v}\right\Vert >\lambda \right\} , \\
S_{3} &=&\inf \left\{ t>0,\left\vert \left\langle M\right\rangle
_{t}\right\vert ^{1/2}>\delta \lambda \right\} ,
\end{eqnarray*}%
with the convention that the infimum of the empty set if $\infty .$ Define
the local martingale $N_{t}=M_{S_{3}\wedge S_{2},\left( t+S_{2}\right)
\wedge S_{3}}$noting that $N_{t}\equiv 0$ on $\left\{ S_{2}=\infty \right\} $%
. It is easy to see that 
\begin{equation}
\sup_{u,v\in \left[ 0,S_{3}\right] }\left\Vert \mathbf{M}_{u,v}\right\Vert
\leq \sup_{u,v\in \left[ 0,S_{3}\wedge S_{2}\right] }\left\Vert \mathbf{M}%
_{u,v}\right\Vert +\sup_{u,v\geq 0}\left\Vert \mathbf{N}_{u,v}\right\Vert .
\end{equation}%
By definition of the relevant stopping times, 
\begin{equation*}
\mathbb{P}\left( \sup_{u,v\geq 0}\left\Vert \mathbf{M}_{u,v}\right\Vert
>\beta \lambda ,\left\vert \left\langle M\right\rangle _{\infty }\right\vert
^{1/2}\leq \delta \lambda \right) =\mathbb{P}\left( S_{1}<\infty
,S_{3}=\infty \right) .
\end{equation*}%
On the event $\left\{ S_{1}<\infty ,S_{3}=\infty \right\} $ one has%
\begin{equation*}
\sup_{u,v\in \left[ 0,S_{3}\right] }\left\Vert \mathbf{M}_{u,v}\right\Vert
>\beta \lambda
\end{equation*}%
and, since $S_{2}\leq S_{1}$, one also has $\sup_{u,v\in \left[
0,S_{3}\wedge S_{2}\right] }\left\Vert \mathbf{M}_{u,v}\right\Vert =\lambda
. $ Hence, on $\left\{ S_{1}<\infty ,S_{3}=\infty \right\} ,$%
\begin{equation*}
\sup_{u,v\geq 0}\left\Vert \mathbf{N}_{u,v}\right\Vert \geq \sup_{u,v\in %
\left[ 0,S_{3}\right] }\left\Vert \mathbf{M}_{u,v}\right\Vert -\sup_{u,v\in %
\left[ 0,S_{3}\wedge S_{2}\right] }\left\Vert \mathbf{M}_{u,v}\right\Vert
\geq \left( \beta -1\right) \lambda .
\end{equation*}%
Therefore, using Proposition \ref{TschebTypeBound},%
\begin{eqnarray*}
\mathbb{P}\left( \sup_{u,v\geq 0}\left\Vert \mathbf{M}_{u,v}\right\Vert
>\beta \lambda ,\left\vert \left\langle M\right\rangle _{\infty }\right\vert
^{1/2}\leq \delta \lambda \right) &\leq &\mathbb{P}\left( \sup_{u,v\geq
0}\left\Vert \mathbf{N}_{u,v}\right\Vert \geq \left( \beta -1\right) \lambda
\right) \\
&\leq &\frac{A}{\left( \beta -1\right) ^{2}\lambda ^{2}}\mathbb{E}\left(
\left\vert \left\langle N\right\rangle _{\infty }\right\vert \right) .
\end{eqnarray*}%
From the definition of $N$, for every $t\in \left[ 0,\infty \right] $,%
\begin{equation*}
\left\langle N\right\rangle _{t}=\left\langle M\right\rangle _{S_{3}\wedge
S_{2},\left( t+S_{2}\right) \wedge S_{3}}.
\end{equation*}%
On $\left\{ S_{2}=\infty \right\} $ we have $\left\langle N\right\rangle
_{\infty }=0$ while on $\left\{ S_{2}<\infty \right\} $ we have, from
definition of $S_{3}$,%
\begin{equation*}
\left\vert \left\langle N\right\rangle _{\infty }\right\vert =\left\vert
\left\langle M\right\rangle _{S_{3}\wedge S_{2},S_{3}}\right\vert
=\left\vert \left\langle M\right\rangle _{S_{3}}-\left\langle M\right\rangle
_{S_{3}\wedge S_{2}}\right\vert \leq 2\left\vert \left\langle M\right\rangle
_{S_{3}}\right\vert =2\delta ^{2}\lambda ^{2}.
\end{equation*}%
It follows that%
\begin{equation*}
\mathbb{E}\left( \left\vert \left\langle N\right\rangle _{\infty
}\right\vert \right) \leq 2\delta ^{2}\lambda ^{2}\mathbb{P}\left(
S_{2}<\infty \right) =2\delta ^{2}\lambda ^{2}\mathbb{P}\left( \sup_{u,v\geq
0}\left\Vert \mathbf{M}_{u,v}\right\Vert >\lambda \right)
\end{equation*}%
and we have the estimate%
\begin{equation*}
\mathbb{P}\left( \sup_{u,v\geq 0}\left\Vert \mathbf{M}_{u,v}\right\Vert
>\beta \lambda ,\left\vert \left\langle M\right\rangle _{\infty }\right\vert
^{1/2}\leq \delta \lambda \right) \leq \frac{2A\delta ^{2}}{\left( \beta
-1\right) ^{2}}\mathbb{P}\left( \sup_{u,v\geq 0}\left\Vert \mathbf{M}%
_{u,v}\right\Vert >\lambda \right) .
\end{equation*}%
An application of the good $\lambda $-inequality finishes the proof.
\end{proof}

\section{Path regularity and $p$-variation BDG}

\bigskip Let $p>2$. From \cite{CL} it is known that for every $\mathbf{M}\in 
\mathcal{M}_{0,\text{\textrm{loc}}}^{c}\left( G^{2}\left( \mathbb{R}%
^{d}\right) \right) $ and every $T>0$ 
\begin{equation}
\left\Vert \mathbf{M}\right\Vert _{p\text{-var;}\left[ 0,T\right] }<\infty 
\text{ a.s.}  \label{FinitePvarForEnhMart}
\end{equation}%
Here, we go one step further and provided quantitative bounds for the $p$%
-variation of the enhanced martingale in terms of $\left\langle
M\right\rangle _{T}$. En passant, we give a simplified proof of (\ref%
{FinitePvarForEnhMart}).

\begin{proposition}
Let $\mathbf{M}\in \mathcal{M}_{0,\text{\textrm{loc}}}^{c}\left( G^{2}\left( 
\mathbb{R}^{d}\right) \right) $. Then, for every $T>0,$%
\begin{equation*}
\left\Vert \mathbf{M}\right\Vert _{p\text{-var;}\left[ 0,T\right] }<\infty 
\text{ a.s.}
\end{equation*}
\end{proposition}

\begin{proof}
There exists a sequence of stopping times $\tau _{n}\rightarrow \infty $
a.s. such that $M^{\tau _{n}}$ and $\left\langle M^{\tau _{n}}\right\rangle $
are bounded (for instance, $\tau _{n}=\inf \{t:\left\vert M_{t}\right\vert
>n $ or $\left\vert \left\langle M\right\rangle _{t}\right\vert >n\}$ will
do.) Since%
\begin{equation*}
\mathbb{P}\left( \left\Vert \mathbf{M}\right\Vert _{p\text{-var;}\left[ 0,T%
\right] }\neq \left\Vert \mathbf{M}\right\Vert _{p\text{-var;}\left[
0,T\wedge \tau _{n}\right] }\right) \leq \mathbb{P}\left( \tau _{n}<T\right)
\rightarrow 0\text{ as }n\rightarrow \infty
\end{equation*}%
it suffices to consider the lift of a bounded continuous martingale with
bounded quadratic variation. We can work with the $l^{1}$-norm on $\mathbb{R}%
^{d}$, $\left\vert a\right\vert =\sum_{i=1}^{d}\left\vert a_{i}\right\vert . 
$ The time change $\phi \left( t\right) :=\inf \left\{ s:\left\vert
\left\langle M\right\rangle _{s}\right\vert >t\right\} $ may have jumps but
continuity of $\left\vert \left\langle M\right\rangle \right\vert $ ensures
that $|\left\langle M\right\rangle _{\phi \left( t\right) }|\ =t$. From
definition of $\phi $ and the BDG inequality on the group, both $\left\vert
\left\langle M\right\rangle \right\vert $ and $\mathbf{M}$ are constant on
the intervals $\left[ \phi _{t-},\phi _{t}\right] $. It follows that $%
\mathbf{X}_{t}=\mathbf{M}_{\phi \left( t\right) }$ defines a continuous%
\footnote{%
From Lemma 2, $\mathbf{X}=S_{2}\left( M\circ \phi \right) ,$ the lift of a
continuous local martingale. In particular, this is another way to see
continuity of $\mathbf{X}$.} path from $\left[ 0,\left\vert \left\langle
M\right\rangle _{T}\right\vert \right] $ to $G^{2}\left( \mathbb{R}%
^{d}\right) $ and it is easy to see that%
\begin{equation*}
\left\Vert \mathbf{X}\right\Vert _{p\text{-var},\left[ 0,\left\vert
\left\langle M\right\rangle _{T}\right\vert \right] }=\left\Vert \mathbf{M}%
\right\Vert _{p\text{-var},\left[ 0,T\right] }.
\end{equation*}%
As argued in the beginning of the proof, we may assume that $\left\vert
\left\langle M\right\rangle _{T}\right\vert \leq R$ for some deterministic $%
R $ large enough. Therefore, 
\begin{eqnarray}
\mathbb{P}\left( \left\Vert \mathbf{M}\right\Vert _{p\text{-var},\left[ 0,T%
\right] }>K\right) &=&\mathbb{P}\left( \left\Vert \mathbf{X}\right\Vert _{p%
\text{-var},\left[ 0,\left\vert \left\langle M\right\rangle _{T}\right\vert %
\right] },\left\vert \left\langle M\right\rangle _{T}\right\vert \leq
R\right)  \label{MltX} \\
&\leq &\mathbb{P}\left( \left\Vert \mathbf{X}\right\Vert _{p\text{-var},%
\left[ 0,R\right] }>K\right) .  \notag
\end{eqnarray}%
We go on to show that $\mathbf{X}$ is in fact H\"{o}lder continuous. For $%
0\leq s\leq t\leq R$, we can use the BDG\ inequality on the group, theorem %
\ref{BDGgroupUniform}, to obtain%
\begin{equation*}
\mathbb{E}\left( \left\Vert \mathbf{X}_{s,t}\right\Vert ^{2q}\right) =%
\mathbb{E}\left( \left\Vert \mathbf{M}_{\phi \left( s\right) ,\phi \left(
t\right) }\right\Vert ^{2q}\right) \leq C_{q}\mathbb{E}\left( \left\vert
\left\langle M\right\rangle _{\phi \left( t\right) }-\left\langle
M\right\rangle _{\phi \left( s\right) }\right\vert ^{q}\right) .
\end{equation*}%
Observe that 
\begin{eqnarray*}
\left\vert \left\langle M\right\rangle _{\phi \left( t\right) }-\left\langle
M\right\rangle _{\phi \left( s\right) }\right\vert &=&\sum_{i}\left(
\left\langle M^{i}\right\rangle _{\phi \left( t\right) }-\left\langle
M^{i}\right\rangle _{\phi \left( s\right) }\right) \\
&=&\left\vert \left\langle M\right\rangle _{\phi \left( t\right)
}\right\vert -\left\vert \left\langle M\right\rangle _{\phi \left( s\right)
}\right\vert =t-s.
\end{eqnarray*}%
Thus, for all $q<\infty $ there exists a constant $C_{q}$ s.t.%
\begin{equation*}
\mathbb{E}\left( \left\Vert \mathbf{X}_{s,t}\right\Vert ^{2q}\right) \leq
C_{q}\left\vert t-s\right\vert ^{q}.
\end{equation*}%
Knowing that $\mathbf{X}$ is continuous, we can apply GRR\footnote{%
There is no modification of $\mathbf{X}$ needed.} for paths in $\left(
G^{2}\left( \mathbb{R}^{d}\right) ,d\right) $ to see that $\left\Vert 
\mathbf{X}\right\Vert _{1/p\text{-H\"{o}lder},\left[ 0,R\right] }\in L^{q}$
for all $q\in \lbrack 1,\infty )$ and 
\begin{equation*}
\mathbb{P}\left( \left\Vert \mathbf{X}\right\Vert _{p\text{-var},\left[ 0,R%
\right] }>K\right) \leq \frac{\mathbb{E}\left( \left\Vert \mathbf{X}%
\right\Vert _{p\text{-var},\left[ 0,R\right] }\right) }{K}\leq \frac{\mathbb{%
E}\left( \left\Vert \mathbf{X}\right\Vert _{1/p\text{-H\"{o}lder},\left[ 0,R%
\right] }\right) }{K}
\end{equation*}%
tends to zero as $K\rightarrow \infty $. Together with (\ref{MltX}) we see
that $\left\Vert \mathbf{M}\right\Vert _{p\text{-var},\left[ 0,T\right]
}<\infty $ with probability $1$ as claimed.
\end{proof}

We are now going to prove a $p$-variation version of BDG. For $\mathbb{R}$%
-valued martingales this result is due to L\'{e}pingle, \cite{L}. With the
preparations made, our proof follows the same lines.

\begin{lemma}
There exists a constant $A$ such that for all continous local martingales $M$%
, for all $\lambda >0,$ 
\begin{equation*}
\mathbb{P}\left( \left\Vert \mathbf{M}\right\Vert _{p\text{-var;}[0,\infty
)}>\lambda \right) \leq A\frac{\mathbb{E}\left( \left\vert \left\langle
M\right\rangle _{\infty }\right\vert \right) }{\lambda ^{2}}.
\end{equation*}
\end{lemma}

\begin{proof}
If suffices to prove the statement when $\lambda =1$ (the general case
follows by considering $M/\lambda $ with lift $\delta _{1/\lambda }\mathbf{M}
$). The statement then reduces to 
\begin{equation*}
\exists A:\forall M:\mathbb{P}\left[ \left\Vert \mathbf{M}\right\Vert _{p%
\text{-var;}[0,\infty )}>1\right] \leq A\,\mathbb{E}\left( \left\vert
\left\langle M\right\rangle _{\infty }\right\vert \right) .
\end{equation*}%
Assume this is false. Then for every $A$, and in particular for $A\left(
k\right) \equiv k^{2}$,there exists $M\equiv M^{\left( k\right) }$ with lift 
$\mathbf{M}^{\left( k\right) }$ s.t. the condition is violated,%
\begin{equation*}
k^{2}\,\mathbb{E}\left[ \left\vert \left\langle M^{\left( k\right)
}\right\rangle _{\infty }\right\vert \right] <\mathbb{P}\left[ \left\Vert 
\mathbf{M}^{\left( k\right) }\right\Vert _{p\text{-var;}[0,\infty )}>1\right]
\leq 1.
\end{equation*}%
Set $u_{k}=\mathbb{P}\left[ \left\Vert \mathbf{M}^{\left( k\right)
}\right\Vert _{p\text{-var;}[0,\infty )}>1\right] $, $n_{k}=\left[ 1/u_{k}+1%
\right] \in \mathbb{N}$ and note that $1\leq n_{k}u_{k}\leq 2$. Take $n_{k}$
copies of each $M^{\left( k\right) }$ and get a sequence of martingales of
form%
\begin{equation*}
\left( \tilde{M}\right) \equiv (\underset{n_{1}}{\underbrace{M^{\left(
1\right) },...,M^{\left( 1\right) }}};\underset{n_{2}}{\underbrace{M^{\left(
2\right) },...,M^{\left( 2\right) }}};M^{\left( 3\right) },...).
\end{equation*}%
Then%
\begin{equation*}
n_{k}k^{2}\mathbb{E}\left[ \left\vert \left\langle M^{\left( k\right)
}\right\rangle _{\infty }\right\vert \right] \leq n_{k}\mathbb{P}\left[
\left\Vert \mathbf{M}^{\left( k\right) }\right\Vert _{p\text{-var;}[0,\infty
)}>1\right] =n_{k}u_{k}\leq 2.
\end{equation*}%
and%
\begin{equation*}
\sum_{k}\mathbb{P}\left[ \left\Vert \mathbf{\tilde{M}}^{\left( k\right)
}\right\Vert _{p\text{-var;}[0,\infty )}>1\right] =\sum_{k}n_{k}u_{k}=+\infty
\end{equation*}%
while%
\begin{equation*}
\sum_{k}\mathbb{E}\left[ \left\vert \left\langle \tilde{M}^{\left( k\right)
}\right\rangle _{\infty }\right\vert \right] =\sum_{k}n_{k}\mathbb{E}\left[
\left\vert \left\langle M^{\left( k\right) }\right\rangle _{\infty
}\right\vert \right] \leq \sum_{k}\frac{2}{k^{2}}<\infty .
\end{equation*}%
Thus, if the claimed statement is false, there exists a sequence of
martingales, we now revert to write $M^{\left( k\right) },\mathbf{M}^{\left(
k\right) }$ instead of $\tilde{M}^{\left( k\right) },\mathbf{\tilde{M}}%
^{\left( k\right) }$ respectively, each defined on some filtered probability
space $\left( \Omega ^{k},\left( \mathcal{F}_{t}^{k}\right) ,\mathbb{P}%
^{k}\right) $ with the two properties%
\begin{equation*}
\sum_{k}\mathbb{P}^{k}\left[ \left\Vert \mathbf{M}^{\left( k\right)
}\right\Vert _{p\text{-var;}[0,\infty )}>1\right] =+\infty \text{ and }%
\sum_{k}\mathbb{E}^{k}\left[ \left\vert \left\langle M^{\left( k\right)
}\right\rangle _{\infty }\right\vert \right] <\infty .
\end{equation*}%
Define the probability space $\Omega =\bigotimes_{k=1}^{\infty }\Omega ^{k},$
the probability $\mathbb{P}=\bigotimes_{k=1}^{\infty }\mathbb{P}^{k},$ and
the filtration $\left( \mathcal{F}_{t}\right) $ on $\Omega $ given by%
\begin{equation*}
\mathcal{F}_{t}=\left( \bigotimes_{i=1}^{k-1}\mathcal{F}_{\infty
}^{i}\right) \otimes \mathcal{F}_{g\left( k-t\right) }^{k}\otimes \left(
\bigotimes_{j=k+1}^{\infty }\mathcal{F}_{0}^{k}\right) \text{ \ for }k-1\leq
t<k.
\end{equation*}%
where $g\left( u\right) =1/u-1$ maps $\left[ 0,1\right] \rightarrow \left[
0,\infty \right] $. Then, a continous martingale on $\left( \Omega ,\left( 
\mathcal{F}_{t}\right) ,\mathbb{P}\right) $ is defined by concatenation, \ 
\begin{equation*}
M_{t}=\sum_{i=1}^{k-1}M_{\infty }^{\left( i\right) }+M_{g\left( k-t\right)
}^{\left( k\right) }\,\,\,\text{for \ \ }k-1\leq t<k.
\end{equation*}%
which implies%
\begin{equation*}
\mathbf{M}_{t}=\left( \bigotimes_{i=1}^{k-1}\mathbf{M}_{\infty }^{\left(
i\right) }\right) \otimes \mathbf{M}_{g\left( k-t\right) }^{\left( k\right)
}.
\end{equation*}%
We also observe that, again for $k-1\leq t<k$, 
\begin{equation*}
\left\langle M\right\rangle _{t}=\sum_{i=1}^{k-1}\left\langle M^{\left(
i\right) }\right\rangle _{\infty }+\left\langle M^{\left( k\right)
}\right\rangle _{g\left( k-t\right) }.
\end{equation*}%
In particular, $\left\langle M\right\rangle _{\infty }=\sum_{k}\left\langle
M^{\left( k\right) }\right\rangle _{\infty }$ and, using the second property
of the martingale sequence, $\mathbb{E}\left( \left\vert \left\langle
M\right\rangle _{\infty }\right\vert \right) <\infty $. Define the events%
\begin{equation*}
A_{k}=\left\{ \left\Vert \mathbf{M}\right\Vert _{p\text{-var;}\left[ k-1,k%
\right] }>1\right\} .
\end{equation*}%
Then, using the first property of the martingale sequence, 
\begin{equation*}
\sum_{k}\mathbb{P}\left( A_{k}\right) =\sum_{k}\mathbb{P}^{k}\left(
\left\Vert \mathbf{M}^{k}\right\Vert _{p\text{-var;}[0,\infty )}>1\right)
=\infty .
\end{equation*}%
Since the events $\left\{ A_{k}:k\geq 1\right\} $ are independent, the
Borel-Cantelli lemma implies that $\mathbb{P}\left( A_{k}\text{ i.o.}\right)
=1$. Thus, almost surely, for all $K>0$ there exists a finite number of
increasing times $t_{0},\cdots ,t_{n}\in \lbrack 0,\infty )$ so that%
\begin{equation*}
\sum_{i=1}^{n}\left\Vert \mathbf{M}_{t_{i-1},t_{i}}\right\Vert >K
\end{equation*}%
and $\left\Vert \mathbf{M}\right\Vert _{p\text{-var;}[0,\infty )}$ must be
equal to $+\infty $ with probability one. We now define a martingale $N$ by
time-change, namely via $f\left( t\right) =t/\left( 1-t\right) $ for $0\leq
t<1$ and $f\left( t\right) =\infty $ for $t\geq 1$,%
\begin{equation*}
N:t\mapsto M_{f\left( t\right) }.
\end{equation*}%
Note that $\mathbb{E}\left( \left\vert \left\langle M\right\rangle _{\infty
}\right\vert \right) <\infty $ so that $M$ can be extended to a (continuous)
martingale indexed by $\left[ 0,\infty \right] $ and $N$ is indeed a
continuous martingale with lift $\mathbf{N}$. Since lifts interchange with
time changes, $\left\Vert \mathbf{N}\right\Vert _{p\text{-var;}%
[0,1]}=\left\Vert \mathbf{M}\right\Vert _{p\text{-var;}[0,\infty )}=+\infty $
with probability one. But this contradicts to $p$-variation regularity
result above.
\end{proof}

The very same argument that was used in the proof of Theorem \ref%
{BDGgroupUniform}\ now leads to the following BDG inequality for enhanced
continuous local martingales w.r.t. homogenuous $p$-variation norm.

\begin{theorem}
\label{pVarBDG}Let $F$ be a moderate function, $\mathbf{M}\in \mathcal{M}_{0,%
\text{\textrm{loc}}}^{c}\left( G^{2}\left( \mathbb{R}^{d}\right) \right) ,$
and $\left\vert \cdot \right\vert ,\,\left\Vert \cdot \right\Vert $
continuous homogonous norm on $\mathbb{R}^{d},G^{2}\left( \mathbb{R}%
^{d}\right) $ respectively and $p>2$. Then there exists a constant $%
C=C\left( p,F,d,\left\vert \cdot \right\vert ,\left\Vert \cdot \right\Vert
\right) $ so that%
\begin{equation*}
C^{-1}\mathbb{E}\left( F\left( \left\vert \left\langle M\right\rangle
_{\infty }\right\vert ^{1/2}\right) \right) \leq \mathbb{E}\left( F\left(
\left\Vert \mathbf{M}\right\Vert _{p\text{-var;}[0,\infty )}\right) \right)
\leq C\mathbb{E}\left( F\left( \left\vert \left\langle M\right\rangle
_{\infty }\right\vert ^{1/2}\right) \right) .
\end{equation*}
\end{theorem}

\begin{remark}
When $p\in \left( 2,3\right) $ and $N\in \left\{ 3,4,...\right\} $, $\mathbf{%
M}$ lifts uniquely to a $G^{N}\left( \mathbb{R}^{d}\right) $-valued path
with finite homogenuous $p$-variation regularity, denoted by $S_{N}\left( 
\mathbf{M}\right) $, which is identified with the first $N$ iterated
Stratonovich integrals of $M$. A basic theorem of Lyons asserts that%
\begin{equation*}
\left\Vert S_{N}\left( \mathbf{M}\right) \right\Vert _{p\text{-var}}\leq
C\left( N\right) \left\Vert \mathbf{M}\right\Vert _{p\text{-var}}
\end{equation*}%
and BDG inequalities for the $p$-variation of this step-$N$ lift are an
immediate corollary of Theorem \ref{pVarBDG}.
\end{remark}

\section{Approximations}

We now only consider (lifted) local martingales on $\left[ 0,T\right] $,
defined or identified with local martingales stopped at $T>0.$

\subsection{Geodesic approxiations}

The $p$-variation norm of geodesics approximations is uniformly controlled
by the original $p$-variation norm. Therefore%
\begin{equation*}
\mathbb{E}\left( F\left( \sup_{D}\left\Vert \mathbf{M}^{\left[ D\right]
}\right\Vert _{p\text{-var;}[0,T]}\right) \right) \leq C\mathbb{E}\left(
F\left( \left\vert \left\langle M\right\rangle _{T}\right\vert ^{1/2}\right)
\right)
\end{equation*}%
where $\mathbf{M}^{\left[ D\right] }$ denotes the geodesics approxiation to $%
\mathbf{M}$ based on some dissection $D$ of $\left[ 0,T\right] $. Note that
this is stronger than%
\begin{equation*}
\sup_{D}\mathbb{E}\left( F\left( \left\Vert \mathbf{M}^{\left[ D\right]
}\right\Vert _{p\text{-var;}[0,T]}\right) \right) \leq C\mathbb{E}\left(
F\left( \left\vert \left\langle M\right\rangle _{T}\right\vert ^{1/2}\right)
\right)
\end{equation*}%
which is what we are going to show for piecewise linear approximations.

\subsection{Piecewise linear approximations}

Let $D=\left( t_{i}\right) $ be a subdivision of $\left[ 0,T\right] .$ Given 
$x\in C\left( \left[ 0,T\right] ,\mathbb{R}^{d}\right) $ we define $x^{D}$
to be the piecewise linear approximation of $x$ which coincides with $x$ on $%
D.$ Since $x^{D}$ is of bounded variation, it admits a canonical lift to a $%
G^{2}\left( \mathbb{R}^{d}\right) $-valued path, denoted by $\mathbf{x}^{D}.$
This notation applies path-by-path to $M\in \mathcal{M}_{0,\text{\textrm{loc}%
}}^{c}\left( \mathbb{R}^{d}\right) $, we write $\mathbf{M}^{D}=\mathbf{M}%
^{D}\left( \omega \right) $ for the lifted piecewise linear approximation to 
$M\left( \omega \right) $. The next lemma involves no probabilty.

\begin{lemma}
\label{pVarOfLiftedPWApprox}Set $\mathbf{x}^{D}=S_{2}\left( x^{D}\right) $
where $x^{D}$ is linear between the points of $D$. Then there exists a
constant $C=C=C\left( d,\left\vert \cdot \right\vert ,\left\Vert \cdot
\right\Vert \right) $ such that%
\begin{equation*}
\left\Vert \mathbf{x}^{D}\right\Vert _{p-var;\left[ 0,T\right] }\leq C\left(
\max_{\left( s_{k}\right) \subset D}\sum_{k}\left\Vert \mathbf{x}%
_{s_{k},s_{k+1}}^{D}\right\Vert ^{p}\right) ^{1/p}+C\left\vert x\right\vert
_{p\text{-var};\left[ 0,T\right] }.
\end{equation*}
\end{lemma}

\begin{proof}
First we note that $\left\Vert \mathbf{x}_{s,t}^{D}\right\Vert ^{p}\leq
3^{p-1}\left[ \left\vert x_{s,s^{D}}^{D}\right\vert ^{p}+\left\Vert \mathbf{x%
}_{s^{D},t_{D}}^{D}\right\Vert ^{p}+\left\vert x_{t_{D},t}^{D}\right\vert
^{p}\right] $. Now let $\left( u_{k}\right) $ be a dissection of $\left[ 0,T%
\right] $, unrelated to $D$. Recall that $u^{D}$ resp. $u_{D}$ refers to the
right- resp. left-neighbours of $u$ in $D$. 
\begin{eqnarray*}
\sum_{k}\left\Vert \mathbf{x}_{u_{k},u_{k+1}}^{D}\right\Vert ^{p} &\leq
&3^{p-1}\sum_{k}\left\Vert \mathbf{x}_{u_{k}^{D},u_{k+1,D}}^{D}\right\Vert
^{p}+3^{p-1}\sum_{k}\left[ \left\vert x_{u_{k},u_{k}^{D}}^{D}\right\vert
^{p}+\left\vert x_{u_{k+1,D},u_{k}}^{D}\right\vert ^{p}\right] \\
&\leq &3^{p-1}\left( \max_{\left( s_{k}\right) \subset D}\sum_{k}\left\Vert 
\mathbf{x}_{s_{k},s_{k+1}}^{D}\right\Vert ^{p}\right) +3^{p-1}\left\vert
x^{D}\right\vert _{p\text{-var;}\left[ 0,T\right] } \\
&\leq &3^{p-1}\left( \max_{\left( s_{k}\right) \subset D}\sum_{k}\left\Vert 
\mathbf{x}_{s_{k},s_{k+1}}^{D}\right\Vert ^{p}\right) +C3^{p-1}\left\vert
x\right\vert _{p\text{-var;}\left[ 0,T\right] }.
\end{eqnarray*}
\end{proof}

\begin{theorem}
\label{UpperBDGboundUniformOverDissections}Let $F$ be a moderate function, $%
\mathbf{M}$ $\in \mathcal{M}_{0,\text{\textrm{loc}}}^{c}\left( G^{2}\left( 
\mathbb{R}^{d}\right) \right) ,$ and $\left\vert \cdot \right\vert
,\,\left\Vert \cdot \right\Vert $ continuous homogonous norm on $\mathbb{R}%
^{d},G^{2}\left( \mathbb{R}^{d}\right) $ respectively. Then there exists a
constant $C=C\left( p,F,d,\left\vert \cdot \right\vert ,\left\Vert \cdot
\right\Vert \right) $ so that for all dissections $D$ of $\left[ 0,T\right]
, $%
\begin{equation*}
\mathbb{E}\left( F\left( \left\Vert \mathbf{M}^{D}\right\Vert _{p-var\text{;}%
[0,T]}\right) \right) \leq C\mathbb{E}\left( F\left( \left\vert \left\langle
M\right\rangle _{T}\right\vert ^{1/2}\right) \right) .
\end{equation*}
\end{theorem}

\begin{proof}
From Lemma \ref{pVarOfLiftedPWApprox}, $\left\Vert \mathbf{M}^{D}\right\Vert
_{p\text{-var;}\left[ 0,T\right] }$ is bounded by $C\left\vert M\right\vert
_{p\text{-var;}\left[ 0,T\right] }$ plus%
\begin{eqnarray*}
&&C\left( \max_{\left( s_{k}\right) \subset D}\sum_{k}\left\Vert \mathbf{M}%
_{s_{k},s_{k+1}}^{D}\right\Vert ^{p}\right) ^{1/p}\leq C\left( \max_{\left(
s_{k}\right) \subset D}\sum_{k}\left\Vert \mathbf{M}_{s_{k},s_{k+1}}\right%
\Vert ^{p}\right) ^{1/p} \\
&&+C\left( \max_{\left( s_{k}\right) \subset D}\sum_{k}d\left( \mathbf{M}%
_{s_{k},s_{k+1}},\mathbf{M}_{s_{k},s_{k+1}}^{D}\right) ^{p}\right) ^{1/p}.
\end{eqnarray*}%
Trivially, $\left\vert M\right\vert _{p\text{-var;}\left[ 0,T\right] }\leq
\left\Vert \mathbf{M}\right\Vert _{p\text{-var;}\left[ 0,T\right] }$ and
with a new constant $C$, 
\begin{equation*}
\left\Vert \mathbf{M}^{D}\right\Vert _{p\text{-var;}\left[ 0,T\right] }\leq
C\left\Vert \mathbf{M}\right\Vert _{p\text{-var;}\left[ 0,T\right] }+C\left(
\max_{\left( s_{k}\right) \subset D}\sum_{k}d\left( \mathbf{M}%
_{s_{k},s_{k+1}},\mathbf{M}_{s_{k},s_{k+1}}^{D}\right) ^{p}\right) ^{1/p}%
\text{.}
\end{equation*}%
For fixed $k,$ there are $i<j$ so that $s_{k}=t_{i}$ and $s_{k+1}=t_{j}$.
Then%
\begin{eqnarray*}
\mathbf{M}_{s_{k},s_{k+1}} &=&\bigotimes_{l=i}^{j-1}\exp \left(
M_{t_{l},t_{l+1}}+A_{t_{l},t_{l+1}}\right) \\
\,\mathbf{M}_{s_{k},s_{k+1}}^{D}\, &=&\bigotimes_{l=i}^{j-1}\exp \left(
M_{t_{l},t_{l+1}}\right) .
\end{eqnarray*}%
Hence, $d\left( \mathbf{M}_{s_{k},s_{k+1}},\mathbf{M}_{s_{k},s_{k+1}}^{D}%
\right) $ equals%
\begin{equation}
\left\Vert \mathbf{M}_{s_{k},s_{k+1}}^{-1}\otimes \mathbf{M}%
_{s_{k},s_{k+1}}^{D}\right\Vert =\left\Vert \exp \left(
\sum_{l=i}^{j-1}A_{t_{l},t_{l+1}}\right) \right\Vert \leq C\left\vert
\sum_{l=i}^{j-1}A_{t_{l},t_{l+1}}\right\vert ^{1/2}.  \label{SumOfSmallAreas}
\end{equation}%
The key idea is to introduce the (vector-valued) discrete-time martingale%
\begin{equation*}
Y_{j}=\sum_{l=0}^{j-1}A_{t_{l},t_{l+1}}\in so\left( d\right) .
\end{equation*}%
From (\ref{SumOfSmallAreas}) and equivalence of homogenous norms we have%
\begin{equation*}
\max_{\left( s_{k}\right) \subset D}\sum_{k}d\left( \mathbf{M}%
_{s_{k},s_{k+1}},\mathbf{M}_{s_{k},s_{k+1}}^{D}\right) ^{p}\leq
C\max_{\left\{ i_{1},...,i_{n}\right\} \subset \left\{ 1,...,\#D\right\}
}\sum_{k}\left\vert Y_{i_{k+1}}-Y_{i_{k}}\right\vert ^{p/2},
\end{equation*}%
which leads to%
\begin{eqnarray*}
\left\Vert \mathbf{M}^{D}\right\Vert _{p\text{-var}} &\leq &C\left\Vert 
\mathbf{M}\right\Vert _{p\text{-var}}+C\sqrt{\left( \max_{\left\{
i_{1},...,i_{n}\right\} \subset \left\{ 1,...,\#D\right\}
}\sum_{k}\left\vert Y_{i_{k+1}}-Y_{i_{k}}\right\vert ^{p/2}\right) ^{2/p}} \\
&=&C\left\Vert \mathbf{M}\right\Vert _{p\text{-var}}+C\sqrt{\left\vert
Y\right\vert _{p/2\text{-var}}}.
\end{eqnarray*}%
Using basic properties of moderate functions we have%
\begin{eqnarray*}
\mathbb{E}\left[ F\left( \left\Vert \mathbf{M}^{D}\right\Vert _{p\text{-var}%
}\right) \right] &\leq &C\mathbb{E}\left[ F\left( \left\Vert \mathbf{M}%
\right\Vert _{p\text{-var}}\right) \right] +C\mathbb{E}\left[ F\left( \sqrt{%
\left\vert Y\right\vert _{p/2\text{-var}}}\right) \right] \\
&=&C\mathbb{E}\left[ F\left( \left\Vert \mathbf{M}\right\Vert _{p\text{-var}%
}\right) \right] +C\mathbb{E}\left[ F\circ \sqrt{\cdot }\left( \left\vert
Y\right\vert _{p/2\text{-var}}\right) \right] .
\end{eqnarray*}%
Note that $F\circ \sqrt{\cdot }$ is moderate since $F$ is moderate. Let $%
2<p^{\prime }<p<3$. Then $1<p^{\prime }/2\leq p/2\leq 2$ and and Lemma \ref%
{Prop2bLep76} yields 
\begin{eqnarray*}
\mathbb{E}\left[ F\circ \sqrt{\cdot }\left( \left\vert Y\right\vert _{p/2%
\text{-var}}\right) \right] &\leq &\mathbb{E}\left[ F\circ \sqrt{\cdot }%
\left( \left( \sum_{l}\left\vert Y_{l+1}-Y_{l}\right\vert ^{p^{\prime
}/2}\right) ^{2/p^{\prime }}\right) \right] \\
&=&\mathbb{E}\left[ F\circ \sqrt{\cdot }\left( \left( \sum_{l}\left\vert
A_{t_{l},t_{l+1}}\right\vert ^{p^{\prime }/2}\right) ^{2/p^{\prime }}\right) %
\right] \\
&\leq &\mathbb{E}\left[ F\left( \left( \sum_{l}\left\Vert \mathbf{M}%
_{t_{l},t_{l+1}}\right\Vert ^{p^{\prime }}\right) ^{1/p^{\prime }}\right) %
\right] \\
&\leq &\mathbb{E}\left[ F\left( \left\Vert \mathbf{M}\right\Vert _{p^{\prime
}\text{-var;}\left[ 0,T\right] }\right) \right] .
\end{eqnarray*}%
Combing the last two estimates and using Theorem \ref{pVarBDG} (with $%
p^{\prime }=1+p/2>2$ and $p$ respectively) gives%
\begin{eqnarray*}
\mathbb{E}\left[ F\left( \left\Vert \mathbf{M}^{D}\right\Vert _{p\text{-var;}%
\left[ 0,T\right] }\right) \right] &\leq &C\mathbb{E}\left[ F\left(
\left\Vert \mathbf{M}\right\Vert _{p\text{-var;}\left[ 0,T\right] }\right) %
\right] +C\mathbb{E}\left[ F\left( \left\Vert \mathbf{M}\right\Vert
_{p^{\prime }\text{-var;}\left[ 0,T\right] }\right) \right] \\
&\leq &2C\mathbb{E}\left( F\left( \left\vert \left\langle M\right\rangle
_{T}\right\vert ^{1/2}\right) \right) .
\end{eqnarray*}
\end{proof}

\begin{remark}
We don't expect a lower BDG bound uniformly over all dissections $D$ of $%
\left[ 0,T\right] $. For instance,%
\begin{equation*}
C^{-1}\mathbb{E}\left( F\left( \left\vert \left\langle M\right\rangle
_{T}\right\vert ^{1/2}\right) \right) \leq \mathbb{E}\left( F\left(
\left\vert M^{D}\right\vert _{\infty \text{;}[0,T]}\right) \right)
\end{equation*}%
can't hold since $D=\left\{ 0,T\right\} $ implies $M_{\infty \text{;}%
[0,T]}^{D}=\left\vert M_{T}\right\vert $ and for $F\left( x\right) =x$ we
would control%
\begin{equation*}
\mathbb{E}\left( \left\vert M\right\vert _{\infty \text{;}[0,T]}\right) \sim 
\mathbb{E}\left( \left\vert \left\langle M\right\rangle
_{T}^{1/2}\right\vert \right)
\end{equation*}%
in terms of $\mathbb{E}\left( \left\vert M_{T}\right\vert \right) $ which is
Doob's $L^{q}$ maximal inequality with $q=1$. But, as is well known, one
needs $q>1$ for Doob's $L^{q}$-inequality to hold true.
\end{remark}

Let us now bound the supremum distance between $\mathbf{M}$ and $\mathbf{M}%
^{D}:$

\begin{lemma}
\label{convergencedinfty}Assume that $M$ is a martingale such that%
\begin{equation}
\left\vert M\right\vert _{\infty ;[0,T]}\in L^{q}\left( \Omega \right) \text{
for some }q\geq 1\text{.}  \label{MinfinityLqcondition}
\end{equation}%
\ \ If $D^{n}$ is a sequence of subdivisions whose time steps tends to $0$
when $n$ tends to $\infty $, then $d_{\infty ;\left[ 0,T\right] }\left( 
\mathbf{M,M}^{D_{n}}\right) $ converges to $0$ in $L^{q}$.
\end{lemma}

\begin{remark}
If $q>1$, Doob's maximal inequality implies that (\ref{MinfinityLqcondition}%
) holds for any $L^{q}$-martingale.
\end{remark}

\begin{proof}[Proof of Lemma \protect\ref{convergencedinfty}]
As in the proof of Theorem \ref{UpperBDGboundUniformOverDissections},
equation (\ref{SumOfSmallAreas}) more specifically, we have that when $%
t=t_{i}\in D$%
\begin{equation*}
d\left( \mathbf{M}_{t},\mathbf{M}_{t}^{D}\right) \leq C\left\vert
\sum_{k=0}^{i-1}A_{t_{k},t_{k+1}}\right\vert ^{1/2}.
\end{equation*}%
Next, consider $t\in \left[ t_{i},t_{i+1}\right] $ for some $i$. The path $%
M_{\cdot }^{D}$ restricted to $\left[ t_{i},t_{i+1}\right] $ is a straight
line with no area, hence

\begin{equation*}
\mathbf{M}_{t_{i},t}^{D}=\exp \left( \frac{t-s}{t_{i+1}-t_{i}}%
M_{t_{i},t_{i+1}}\right) .
\end{equation*}

and%
\begin{eqnarray*}
d\left( \mathbf{M}_{t,},\mathbf{M}_{t}^{D}\right) &=&d\left( \mathbf{M}%
_{t_{i}}\otimes \mathbf{M}_{t_{i},t},\mathbf{M}_{t_{i}}^{D}\otimes \mathbf{M}%
_{t_{i},t}^{D}\right) \\
&=&\left\Vert \left( \mathbf{M}_{t_{i},t}^{D}\right) ^{-1}\otimes \left( 
\mathbf{M}_{t_{i}}^{D}\right) ^{-1}\otimes \mathbf{M}_{t_{i}}\otimes \mathbf{%
M}_{t_{i},t}\right\Vert \\
&\leq &\left\Vert \left( \mathbf{M}_{t_{i},t}^{D}\right) \right\Vert
+\left\Vert \left( \mathbf{M}_{t_{i}}^{D}\right) ^{-1}\otimes \mathbf{M}%
_{t_{i}}\right\Vert +\left\Vert \mathbf{M}_{t_{i},t}\right\Vert \\
&\leq &2\sup_{u,v\in \left[ t_{i},t_{i+1}\right] }\left\Vert \mathbf{M}%
_{u,v}\right\Vert +C\max_{i,j}\left\vert
\sum_{l=i}^{j-1}A_{t_{l},t_{l+1}}\right\vert ^{1/2}.
\end{eqnarray*}

For the $L^{q}$ convergence, because $\mathbf{M}$ is almost surely
continuous (in fact, uniformly continuous on the compact $\left[ 0,T\right] $%
)%
\begin{equation*}
\max_{i=0,...,\#D-1}\sup_{s,t\in \left[ t_{i}^{n},t_{i+1}^{n}\right]
}\left\Vert \mathbf{M}_{t_{i},t_{i+1}}\right\Vert \rightarrow 0\text{ a.s.}
\end{equation*}%
Hence, by dominated convergence, 
\begin{equation*}
\lim_{\left\vert D_{n}\right\vert \rightarrow 0}\mathbb{E}\left(
\max_{i=0,...,\#D-1}\sup_{s,t\in \left[ t_{i}^{n},t_{i+1}^{n}\right]
}\left\Vert \mathbf{M}_{t_{i},t_{i+1}}\right\Vert ^{q}\right) =0.
\end{equation*}%
With $Y$ defined as in the proof of Theorem \ref%
{UpperBDGboundUniformOverDissections}, 
\begin{equation*}
\max_{i,j}\left\vert \sum_{l=i}^{j-1}A_{t_{l},t_{l+1}}\right\vert ^{1/2}\leq
C\left[ \left( \max_{\left\{ i_{1},...,i_{n}\right\} \subset \left\{
1,...,\#D\right\} }\sum_{k}\left\vert Y_{i_{k+1}}-Y_{i_{k}}\right\vert
^{p/2}\right) ^{2/p}\right] ^{1/2}
\end{equation*}%
the computation given therein with $F\left( x\right) =x^{q}$ shows%
\begin{eqnarray*}
\mathbb{E}\left( \max_{i,j}\left\vert
\sum_{l=i}^{j-1}A_{t_{l},t_{l+1}}\right\vert ^{q/2}\right) &\leq &C\mathbb{E}%
\left[ F\circ \sqrt{}\left( \max_{\left\{ i_{1},...,i_{n}\right\} \subset
\left\{ 1,...,\#D\right\} }\sum_{k}\left\vert
Y_{i_{k+1}}-Y_{i_{k}}\right\vert ^{p/2}\right) ^{2/p}\right] \\
&\leq &C\mathbb{E}\left[ F\left( \left( \sum_{l:t_{l}\in D_{n}}\left\Vert 
\mathbf{M}_{t_{l},t_{l+1}}\right\Vert ^{q}\right) ^{1/q}\right) \right] \\
&=&\mathbb{E}\left[ \left( \sum_{l:t_{l}\in D_{n}}\left\Vert \mathbf{M}%
_{t_{l},t_{l+1}}\right\Vert ^{q}\right) \right] .
\end{eqnarray*}%
Bu this last expression tends to zero, combining the bounded convergence
theorem with a.s. convergence%
\begin{equation*}
\lim_{n\rightarrow \infty }\sum_{l:t_{l}\in D_{n}}\left\Vert \mathbf{M}%
_{t_{l},t_{l+1}}\right\Vert ^{q}=0.
\end{equation*}%
Indeed, this follows from $\mathbf{M}\in C^{0,q\text{-var}}$ since $q>2$ and
using the usual squeezing argument. To show $L^{q\text{ }}$convergence with
respect to $d_{\infty \text{ }}=d_{\infty ;\left[ 0,T\right] }$, we also
write $\left\Vert \mathbf{\cdot }\right\Vert _{\infty }=\left\Vert \mathbf{%
\cdot }\right\Vert _{\infty ;\left[ 0,T\right] }$ here, recall that 
\begin{equation*}
d_{\infty }\left( \mathbf{M,M}^{D}\right) \leq \sup_{t\in \left[ 0,T\right]
}d\left( \mathbf{M}_{t}\mathbf{,M}_{t}^{D}\right) +c\left\vert \left\Vert 
\mathbf{M}\right\Vert _{\infty }\sup_{t\in \left[ 0,T\right] }d\left( 
\mathbf{M}_{t}\mathbf{,M}_{t}^{D}\right) \right\vert ^{1/2}.
\end{equation*}%
We just showed that $\sup_{t\in \left[ 0,T\right] }d\left( \mathbf{M}_{t}%
\mathbf{,M}_{t}^{D}\right) \rightarrow 0$ in $L^{q}$. Then%
\begin{eqnarray*}
&&\mathbb{E}\left( \left\vert \left\Vert \mathbf{M}\right\Vert _{\infty
}\sup_{t\in \left[ 0,T\right] }d\left( \mathbf{M}_{t}\mathbf{,M}%
_{t}^{D}\right) \right\vert ^{q/2}\right) \\
&\leq &\left( \mathbb{E}\left( \left\vert \left\Vert \mathbf{M}\right\Vert
_{\infty }\right\vert ^{q}\right) \right) ^{1/2}\left( \mathbb{E}\left(
\left\vert \sup_{t\in \left[ 0,T\right] }d\left( \mathbf{M}_{t}\mathbf{,M}%
_{t}^{D}\right) \right\vert ^{q}\right) \right) ^{1/2}.
\end{eqnarray*}%
(Note that by the our BDG\ inqualities%
\begin{equation*}
\mathbb{E}\left( \left\vert \left\Vert \mathbf{M}\right\Vert _{\infty ;\left[
0,T\right] }\right\vert ^{q}\right) \leq C\mathbb{E}\left( \left\vert
\left\langle M\right\rangle _{T}\right\vert ^{q}\right) \leq C\mathbb{E}%
\left( \left\vert \left\vert M\right\vert _{\infty ;\left[ 0,T\right]
}\right\vert ^{q}\right)
\end{equation*}%
and the last expression is finite by assumption.)
\end{proof}

\begin{theorem}
Let $M$ be as in Lemma \ref{convergencedinfty}. Then, $d_{p\text{-var;}\left[
0,T\right] }\left( \mathbf{M}^{D},\mathbf{M}\right) $ converges to $0$ in $%
L^{q}.$ If $M$ is a local martingale, then convergence holds in probability.
\end{theorem}

\begin{proof}
The result for the local martingale will hold if the first result holds, by
a localisation argument that we leave to the reader. We already saw that $%
L^{q}$-convergence holds w.r.t. $d_{\infty }=d_{\infty ;\left[ 0,T\right] }$%
. To go further, writing $d_{p\text{-var}}\equiv d_{p\text{-var;}\left[ 0,T%
\right] }$, we use the interpolation formula 
\begin{equation*}
d_{p\text{-var}}\left( \mathbf{M},\mathbf{M}^{D}\right) \leq Cd_{\infty
}\left( \mathbf{M,M}^{D}\right) ^{1-\frac{p^{\prime }}{p}}\left( \left\Vert 
\mathbf{M}\right\Vert _{p^{\prime }-var}^{\frac{p^{\prime }}{p}}+\left\Vert 
\mathbf{M}^{D}\right\Vert _{p^{\prime }-var}^{\frac{p^{\prime }}{p}}\right)
,\,\,\,2<p^{\prime }<p.
\end{equation*}%
Hence,%
\begin{equation*}
\mathbb{E}\left( \left\vert d_{p-var}\left( \mathbf{M}^{D},\mathbf{M}\right)
\right\vert ^{q}\right) \leq C\mathbb{E}\left( \left( \left\Vert \mathbf{M}%
\right\Vert _{p^{\prime }-var}^{q\frac{p^{\prime }}{p}}+\left\Vert \mathbf{M}%
^{D}\right\Vert _{p^{\prime }-var}^{q\frac{p^{\prime }}{p}}\right) d_{\infty
}\left( \mathbf{M,M}^{D}\right) ^{q\left( 1-\frac{p^{\prime }}{p}\right)
}\right)
\end{equation*}%
Using H\"{o}lder with conjugate exponents $1/\left( p^{\prime }/p\right) $
and $1/\left( 1-p^{\prime }/p\right) $ gives%
\begin{equation*}
\mathbb{E}\left( \left\vert d_{p-var}\left( \mathbf{M}^{D},\mathbf{M}\right)
\right\vert ^{q}\right) \leq C\mathbb{E}\left( \left\Vert \mathbf{M}%
\right\Vert _{p^{\prime }-var}^{q}+\left\Vert \mathbf{M}^{D}\right\Vert
_{p^{\prime }-var}^{q}\right) ^{p^{\prime }/p}\left[ \mathbb{E}\left(
d_{\infty }\left( \mathbf{M,M}^{D}\right) ^{q}\right) \right] ^{1-p^{\prime
}/p}.
\end{equation*}%
But now it suffices to remark, using our BDG estimates, that%
\begin{equation*}
\mathbb{E}\left( \left\Vert \mathbf{M}\right\Vert _{p^{\prime }-var;\left[
0,T\right] }^{q}\right) ,\mathbb{E}\left( \left\Vert \mathbf{M}%
^{D}\right\Vert _{p^{\prime }-var;\left[ 0,T\right] }^{q}\right) \leq C%
\mathbb{E}\left( \left\vert \left\langle M\right\rangle _{T}\right\vert
^{q/2}\right) \leq C\mathbb{E}\left( \left\vert \left\vert M\right\vert
_{\infty ;\left[ 0,T\right] }\right\vert ^{q}\right)
\end{equation*}%
and the last term is finite by assumption.
\end{proof}

\end{document}